\documentclass[11pt,a4paper,reqno,centertags]{article}
\usepackage{xypic,a4wide,amsmath,amsthm,amssymb}

\newtheorem{thm}{Theorem}[section]

\newtheorem{lemma}[thm]{Lemma}
\newtheorem{prop}[thm]{Proposition}

\theoremstyle{remark}
\newtheorem{rem}[thm]{Remark}
\newtheorem{ex}[thm]{Example}

\newenvironment{example}{\begin{ex}\rm}{\qee\end{ex}}
\newcommand{\qee}{\mbox{\hspace{0.2mm}}\hfill$\triangle$}

{\nonumber}
\newcommand{\di}{{\mathrm d}}

\newcommand{\C}{{\mathbb C}}

\begin{document}
\begin{center}
 {\LARGE\bf Gromov-Witten theory of orbicurves, the space of tri-polynomials and Symplectic Field Theory of Seifert fibrations.}\\[15pt]
 {\sc Paolo Rossi}\\
 {(SISSA - Trieste)}
\end{center}

\vspace{1cm}

\begin{abstract}
We compute, with Symplectic Field Theory techniques, the Gromov-Witten theory of $\mathbb{P}^1_{\alpha_1,\ldots,\alpha_a}$, i.e. the complex projective line with $a$ orbifold points. A natural subclass of these orbifolds, the ones with polynomial quantum cohomology, gives rise to a family of (polynomial) Frobenius manifolds and integrable systems of Hamiltonian PDEs, which extend the (dispersionless) bigraded Toda hierarchy (\cite{C}). We then define a  Frobenius structure on the spaces of polynomials in three complex variables of the form $F(x,y,z)=-xyz+P_1(x)+P_2(y)+P_3(z)$ which contains as special cases the ones constructed on the space of Laurent polynomials (\cite{D},\cite{MT}).  We prove a mirror theorem stating that these Frobenius structures are isomorphic to the ones found before for polynomial $\mathbb{P}^1$-orbifolds. Finally we link rational Symplectic Field Theory of Seifert fibrations over $\mathbb{P}^1_{a,b,c}$ with orbifold Gromov-Witten invariants of the base, extending a known result (\cite{B}) valid in the smooth case.
\end{abstract}

\section*{Introduction}

Since the fundamental work of Witten and Kontsevich (\cite{Ko},\cite{W}), Gromov's powerful tool of holomorphic curves and the subsequent successful Gromov-Witten theory showed a remarkable (and quite involved) link between the topology of symplectic manifolds and integrable systems of Hamiltonian PDEs. This phenomenon, studied and encoded into the geometric structure of Frobenius manifolds by Dubrovin (see e.g. \cite{DZ}), provided a bridge to fruitfully exchange insight and results between the two disciplines.

Thanks to the program of Symplectic Field Theory (SFT) initiated by Eliashberg, Givental and Hofer (\cite{E},\cite{EGH}), a new, very suggestive and more direct link has been discovered in the study of holomorphic curves in directed symplectic cobordisms between contact manifolds (or, more in general, manifolds with a stable Hamiltonian structure). The hope is that this theory actually helps in explaining and understanding the above relation. In Symplectic Field Theory the integrable structure encodes some geometric properties at the level of top dimensional stratum of the boundary of the relevant moduli space of curves. This stratum (as in Floer theory) has codimension one, so the Gromov-Witten-like invariants which one considers are directly sensitive to its effects. In the ordinary Gromov-Witten case, on the contrary, these phenomena are hidden in higher codimension, and their effects are much more subtle.

This is a general idea whose details are far from being completely understood. However there are some geometric situations where the common geometric origin of integrable systems from Gromov-Witten and Symplectic Field Theory is quite evident. In particular, in the case of a prequantization bundle $V$ over an integral symplectic manifold $M$ (or, more in general, a Hamiltonian structure of fibration type), as showed in \cite{EGH} and \cite{B}, we have an (at least at genus $0$) explicit correspondence between the integrable hierarchies associated to $V$ (via SFT) and $M$ (via the Frobenius structure of quantum cohomology).

\vspace{1cm}

In this paper we give a procedure for computing the genus $g$ orbifold Gromov-Witten potential of an orbicurve $C$, when only the even (ignoring degree shifting) orbifold cohomology is considered, in terms of certain Hurwitz numbers. In general this potential is a power series of the orbifold cohomology variables $t_0,\ldots,t_l,\mathrm{e}^{t_{l+1}}$ (here $l$ is the dimension of the even orbifold cohomology of the orbicurve). The computation will use a very mild extension of Symplectic Field Theory where we allow orbifold singularities of the target cobordisms. Of course such a extension, in its full generality, would need a thorough study of the moduli space geometry it gives rise, in order to fix the foundations, which is still to be completed even in the ordinary smooth case. However here we consider just the case of $2$-dimensional target cobordisms with $\mathbb{Z}_k$ isolated singularities (in fact punctured orbicurves), where the slight modifications to the theory (mainly the index formula for the dimension of the moduli space of maps and the grading of the variables) and its main results can be guessed easily.

Once we have described the procedure to obtain this ``partial'' potential for $C$ as a series, we restrict to the case of $\mathbb{P}^1_{\alpha_1,\ldots,\alpha_a}$ (the complex projective line with $a$ orbifold points with singularities $\mathbb{Z}_{\alpha_i}$) and we ask when this series truncates. The answer turns out to be that the Frobenius structure is polynomial (in $t_1,\ldots,t_l,\mathrm{e}^{t_{l+1}}$) exactly for any $\mathbb{P}^1_{\alpha_1,\alpha_2}$ (see \cite{MT}) and for $\mathbb{P}^1_{2,2,l-2}$ ($l\geq 4$), $\mathbb{P}^1_{2,3,3}$, $\mathbb{P}^1_{2,3,4}$ and $\mathbb{P}^1_{2,3,5}$. 

Next we construct a family of Frobenius manifolds on the space $M_{p,q,r}$ of polynomials in three complex variables of the form $$F(x,y,z)=-xyz+P_1(x)+P_2(y)+P_3(z)$$
with given degrees $p,q,r$ in $x,y,z$ respectively and $\frac{1}{p}+\frac{1}{q}+\frac{1}{r}>1$. The cases $M_{p,q,1}$ are easily seen to coincide with the spaces of Laurent polynomials of \cite{DZ1}, which have been shown (\cite{MT}) to be isomorphic as Frobenius manifolds to $QH^*_\mathrm{orb}(\mathbb{P}^1_{k,l-k+1})$. This motivates our mirror theorem

\begin{thm}\label{mirror}
Let $M_{p,q,r}$ be the space of tri-polynomials. Then we have the following isomorphisms of Frobenius manifolds:
\begin{equation*}
M_{p,q,r}\cong QH^*_\mathrm{orb}(\mathbb{P}^1_{p,q,r})
\end{equation*}
\end{thm}

\vspace{0.5cm}

An Homological Mirror Symmetry version of this result was conjectured by Takahashi in \cite{T}, where tri-polynomials also appear. Very recently prof. Takahashi communicated he found a proof of this Homological Mirror Symmetry result which should hold for general $p,q,r$. Equivalence of Homological Mirror Symmetry with isomorphism of corresponding Frobenius manifolds is generally believed to be true, but far from being proven. This one can be seen as an instance where both the versions of Mirror Symmetry hold.

\vspace{0.5cm}

Notice, moreover, that $M_{2,2,r}$ can in turn be proven to be isomorphic to $M(D_{r+2},r)$, the Frobenius manifold associated to the extended affine Weyl group $\tilde{D}_l$ (see \cite{DZ1}). This can be promptly done by using Dubrovin and Zhang's reconstruction Theorem $2.1$ of \cite{DZ1}, which states that there is only one polynomial Frobenius structure with the unity vector field, Euler vector field and intersection pairing of $M(D_{r+2},r)$. It is then sufficient to check that $M_{p,q,r}$ has precisely the same $e$, $E$ and intersection metric, the last one given by 
$$(w_1,w_2)^\sim = i_E (w_1 \bullet w_2)$$
where we identify $1$-forms and vectors using the other metric of the Frobenius manifold.

In fact polynomiality follows from the isomorphism $M_{2,2,r}\simeq QH^*_{\mathrm{orb}}(\mathbb{P}^1_{2,2,r})$ and our classification theorem of $\mathbb{P}^1$-orbifolds.

The exceptional cases $E_6$, $E_7$, $E_8$ are even more easily dealt with explicitly computing and comparing the Frobenius potentials. Hence, let $M(R,k)$ be the Frobenius manifold associated to the extended affine Weyl group of the root system $R$ with the choice of the $k$-th root, according to \cite{DZ1}. Then we have the following isomorphisms of Frobenius manifolds:
\begin{equation*}
\begin{split}
M(A_l,k)&\cong QH^*_\mathrm{orb}(\mathbb{P}^1_{k,l-k+1})\qquad \mathrm{see\ \cite{MT}}\\
M(D_l,l-2)&\cong QH^*_\mathrm{orb}(\mathbb{P}^1_{2,2,l-2})\\
M(E_l,4)&\cong QH^*_\mathrm{orb}(\mathbb{P}^1_{2,3,l-3})\qquad l=6,7,8
\end{split}
\end{equation*}

\vspace{0.5cm}

Finally, after extending the result about SFT of Hamiltonian structures of fibration type to the case of an orbifold base $M$, we compute the SFT-potential of a general Seifert fibration over $\mathbb{P}^1_{a,b,c}$ with different contact or Hamiltonian strucures given in the usual way by the fibration itself. This potential will involve, as expected, the integrable systems associated to the rational Gromov-Witten theory of the base. This process can actually be carried on, along the lines of section $2.9.3$ of \cite{EGH}, to extract symplectic and contact invariants of higher dimensional manifolds.

\subsection*{Acknowledgments} 
I am very grateful to my advisor B. Dubrovin for his constant support and expert guidance.

\vspace{0.5cm}

This work is partially supported by the European Science Foundation Programme ``Methods of Integrable Systems, Geometry, Applied Mathematics" (MISGAM), the Marie Curie RTN ``European Network in Geometry, Mathematical Physics and Applications"  (ENIGMA),  and by the Italian Ministry of Universities and Researches (MIUR) research grant PRIN 2006 ``Geometric methods in the theory of nonlinear waves and their applications".

\section{Symplectic Field Theory of punctured Riemann surfaces}

First recall that a Hamiltonian structure (see also \cite{E}) is a pair $(V,\Omega)$, where $V$ is an oriented manifold of dimension $2n-1$ and $\Omega$ a closed $2$-form of maximal rank $2n-2$. The line field $\mathrm{Ker} \Omega$ is called the characteristic line field and we will call characteristic any vector field which generates $\mathrm{Ker} \Omega$. A Hamiltonian structure is called stable if and only if there exists a $1$-form $\lambda$ and a characteristic vector field $R$ (called Reeb) such that $$\lambda(R)=1\qquad\mathrm{and}\qquad i_R \di \lambda=0.$$
A framing of a stable Hamiltonian structure is a pair $(\lambda,J)$ with $\lambda$ as above and $J$ an complex structure on the bundle $\xi=\{\lambda=0\}$ compatible with $\Omega$.

The two main examples of framed Hamiltonian structures we are going to consider arise from contact manifolds and $S^1$-bundles over symplectic manifolds. In the first case, given a contact manifold $(V,\xi=\{\lambda=0\})$ with a compatible complex structure $J$ on $\xi$ we can consider the framed Hamiltonian structure $(V,\Omega=\di \lambda,\lambda,J)$. In the second case, let $(M,\omega)$ be a symplectic manifold with a compatible almost complex structure $J_M$, $p:V\to M$ any $S^1$-bundle and $\lambda$ any $S^1$-connection form over it; then $(V,\Omega=p^*\omega,\lambda,J)$, with $J$ the lift of $J_M$ to the horizontal distribution, is a framed Hamiltonian structure.

\vspace{0.5cm}

We will not recall here the general definition and tools of Symplectic Field Theory. The reader can refer to the preliminary section of \cite{R}, whose notation we conform with, or directly to \cite{EGH}, of which \cite{R} is a quick review.

In this section we prove a general formula for the SFT-potential of a two-dimensional symplectic cobordism $S$ obtained from a genus $g'$ Riemann surface $\Sigma_{g'}$ by removing small discs around $a$ points $z_1,\ldots,z_a\in\Sigma_{g'}$ (see also \cite{R}). We will consider it as a framed symplectic cobordims between the disjoint union of $m$ copies of $S^1$ and the empty set and we will consider only the dependence on the even cohomology classes (i.e those represented by $0$-forms and $2$-forms). The potential will be expressed in terms of Hurwitz numbers $H^{\Sigma_{g'}}_{g,d}(\mu_1,\ldots,\mu_a)$, i.e. the number coverings of $\Sigma_{g'}$ of genus $g$ and degree $d$, branched only over $z_1,\ldots,z_a$ with ramification profile $\mu_1,\ldots,\mu_a$ respectively.
\begin{lemma}\label{multipants}
\begin{equation*}
\begin{split}
\mathbf{F}_S(t_0,s_2,p^1,\ldots,p^a)=&\frac{1}{\hbar}\frac{t_0^2s_2}{2}-\frac{s_2}{24}\\
&+\sum_{g=0}^{\infty}\sum_{d=0}^{\infty}\sum_{|\mu_1|,\ldots,|\mu_a|=d}\mathrm{e}^{d s_2}H^{\Sigma_{g'}}_{g,d}(\mu_1,\ldots,\mu_a)(p^1)^{\mu_1}\ldots(p^a)^{\mu_a}\hbar^{g-1}
\end{split}
\end{equation*}
where $t_0,s_2$ are the variables associated to $\Delta_0=1\in H^0(S)$ and $\Theta_2=\omega\in H_{\mathrm{comp}}^2(S)$ respectively.
\end{lemma}
\begin{proof}
Recall the index theorem for the (virtual) dimension of $\mathcal{M}^A_{g,r,s^+,s^-}$, whence:
$$\text{dim} \mathcal{M}^A_{g,r,s}(S)=2\left(\sum_{i=1}^a s_i+2g-2\right)+2r+2d(2-2g'-a)$$
where $s_i$ is the number of punctures asymptotically mapped to $z_i$ and $s_i=d-\sum_{k=1}^{s_i}(e_k^i-1)$ where $e_k^i$ is the multiplicity of the $k$-th puncture mapped to $z_i$. Let us write the potential for $t_0=s_2=0$: then the relevant component of the moduli space must be zero-dimensional. The above formula, together with Riemann-Hurwitz theorem for branched coverings, ensures then that we have to count curves with no marked points ($r=0$) and no branch values other then $z_1,\ldots,z_n$, hence
$$\mathbf{F}_S(t_0=0,s_2=0,p^1,\ldots,p^a)= \sum_{g=0}^{\infty}\sum_{d=0}^{\infty}\sum_{|\mu_1|,\ldots,|\mu_a|=d}H^{\Sigma_{g'}}_{g,d}(\mu_1,\ldots,\mu_a) (p^1)^{\mu_1}\ldots(p^a)^{\mu_a}\hbar^{g-1}$$
It now remains to apply Theorem $2.7.1$ of \cite{EGH}, making $\mathrm{exp}(\mathbf{F}_S(t_0=0,s_2=0,p^1,\ldots,p^a))$ evolve through Schr\"odinger equation with KdV Hamiltonian $\mathbf{H}=\frac{t_0^2}{2}-\frac{\hbar}{24}+\sum p^i_kq^i_k$ to get the desired formula.
\end{proof}

\section{Gromov-Witten invariants of orbicurves}

We now plan to use Lemma \ref{multipants} and the gluing theorem for composition of cobordisms (Theorem $2.5.3$ of \cite{EGH}) to obtain the orbifold Gromov-Witten potential of a genus $g'$ Riemann surface with orbifold points. In order to do that we will need to consider the mild generalization of Symplectic Field Theory we referred to in the Introduction. Namely we allow the target $2$-dimensional symplectic cobordisms $C$ to have a finite number of isolated codimension-$2$ $\mathbb{Z}_k$ singularities (see e.g. \cite{CR}\cite{CR1} for a review of orbifold geometry). The cobordism $C$ can then be identified with a punctured orbicurve. We will denote by $n$ the number of punctures, by $a$ the number of singular points and by $\alpha_1,\ldots,\alpha_a$ the orders of the singularities of $C$. The relevant moduli space of holomorphic maps will inherit both the characteristics of Chen and Ruan's space of orbicurves (\cite{CR1}) and of the usual SFT's space of curves with punctures asymptotic to Reeb orbits at $\pm \infty$ (\cite{EGH}). We will denote by $\mathcal{M}^d_{g,r,s^+,s^-}(C,J,\mathbf{x})$ the moduli space of holomorphic genus $g$ orbicurves in $C$ of degree $d$ with $r$ marked points, $s^\pm$ positive/negative punctures, which are of type $\mathbf{x}$ in the sense of \cite{CR1} (here $\mathbf{x}$ is a connected component in the inertia orbifold $\tilde{C}$). In particular, smoothness of the cylindrical ends of the target cobordism ensures that the main operations on $\mathcal{M}^d_{g,r,s^+,s^-}(C,J,\mathbf{x})$, first of all its compactification, can be performed completely analogously to the smooth case (with appearance of $n$-story stable curves and the usual rich structure of the boundary).

Let us write, for instance and future reference, the index formula for the dimension of this generalized moduli space:
$$\text{dim}\, \mathcal{M}^d_{g,r,s^+,s^-}(C,J,\mathbf{x})=2(2g+s^+ +s^--2)+2d\,c_1^{\mathrm{rel}}(C)+2r-2\iota(\mathbf{x})$$
where $$c_1^{\mathrm{rel}}(C)=(2-2g-n+\sum_1^a\frac{1-\alpha_i}{\alpha_i})$$ is the first Chern class of $C$, relative to the boundary, and $\iota(\mathbf{x})$ is the degree shift of the connected component $\mathbf{x}$ of $\tilde{C}$ (see \cite{CR} for details).

Similarly, when one considers the Weyl and Poisson graded algebras of SFT (see \cite{EGH} for details) where the Gromov-Witten-like potentials are defined, the grading for the $p$ and $q$ variables, together with the cohomological ones (t=$\sum t_i\theta_i$, where $\theta_i$ is an homogeneous element in $\Omega^*(\tilde{C})$), will be given by formulas that are totally analogous to the ones for the smooth case (see \cite{EGH})
$$\mathrm{deg}\,p_\gamma=+CZ(\gamma)-2$$
$$\mathrm{deg}\,q_\gamma=-CZ(\gamma)-2$$
$$\mathrm{deg}\,t_i=\mathrm{deg}_\mathrm{orb}\,\theta_i-2$$
$$\mathrm{deg}\,z=-2c_1(C)$$
but take into account the effect of the singularities on the orbifold cohomological degree and on the trivialization of the relevant symplectic bundles, giving rise to rational Conley-Zehnder indices. We will give below explicit formulas for the specific cases we are going to use. With this grading, the property of the SFT-potential of a cobordism of being homogeneous of degree $0$, is preserved even in this singular case.

\vspace{0.5cm}

The basic building blocks of our construction for the Gromov-Witten potential of an orbicurve are the potentials for the punctured Riemann surface (Lemma \ref{multipants}) and the \hbox{$\alpha$-orbifold} cap $\mathbb{C}/\mathbb{Z}_\alpha$. Let us compute this last missing element, which we see as a singular symplectic cobordism between $S^1$ and the empty manifold. Consider orbifold cohomology classes of the orbifold cap, which we denote $t_0,t_1,\ldots,t_{\alpha-1}$, where $t_i\in H^{2i/\alpha}_\mathrm{orb}(\mathbb{C}/\mathbb{Z}_\alpha)$. The grading for these variables is given, accordingly, by $\mathrm{deg}\,t_i=\frac{2i}{\alpha}-2$. Moreover, the natural trivialization of the tangent bundle of the cap along the circle boundary gives a fractional CZ-index of $1/\alpha$ for the simple orbit (one actually gets just a natural trivialization for the $\alpha$-fold covering of the cap which gives $CZ(\alpha S^1)=1$), which translates to $\mathrm{deg}\, p_k=-\frac{2k}{\alpha}-2$.
In order to compute the explicit form of the SFT-potential of $\mathbb{C}/\mathbb{Z}_\alpha$, which we denote $\mathbf{F}_\alpha(t_0,\ldots,t_{\alpha-1},p,\hbar)$, we proceed the following way. Recall that $\mathbf{F}_\alpha=\sum\hbar^{g-1}\mathbf{F}_{\alpha,g}$ must be homogeneous of degree $0$, with $\mathrm{deg}\,\hbar=-4$. Since all our variables have negative degree we deduce that the only nonzero terms appear for genus $g=0$. Here we have $$\mathbf{F}_\alpha=\frac{1}{\hbar}\sum_{\substack{i_0,\ldots,i_{\alpha-1}\\ \sum(\alpha-k)i_k=2\alpha}} A_{i_0\ldots i_{\alpha-1}}t_0^{i_0}\ldots t_{\alpha-1}^{i_{\alpha-1}}+\frac{1}{\hbar}\sum_{j=1}^\alpha \sum_{\substack{i_1,\ldots,i_{\alpha-1}\\ \sum(\alpha-k)i_k=\alpha-j}}B_{i_1\ldots i_{\alpha-1};j}\ t_1^{i_1}\ldots t_{\alpha-1}^{i_{\alpha-1}}\,\frac{p_j}{j}$$
One way of determining the coefficients $A$ and $B$ is by gluing any two orbifold caps along their boundary $S^1$ and use Theorems $2.7.1$ and $2.5.3$ of \cite{EGH} to compute the Gromov-Witten potential $\mathbf{F}_{\alpha_1,\alpha_2}(t_0,t_{1,1},\ldots,t_{1,\alpha_1-1},t_{2,1},\ldots,t_{2,\alpha_2-1},s)$ of the resulting $\mathbb{P}^1_{\alpha_1,\alpha_2}$ (here, of course, $s\in H^2(\mathbb{P}^1_{\alpha_1,\alpha_2})$). This potential was already computed by Milanov and Tseng in \cite{MT}, who showed its relation with Carlet's extended bigraded Toda hierarchy \cite{C} and, hence, with extended affine Weyl groups for the root systems $A_l$ (see \cite{DZ1}). Alternatively, one can impose WDVV equations on these potentials and determine the desired coefficients. (this is, by the way, a very efficient way to compute GW-invariants of $\mathbb{P}^1_{\alpha_1,\alpha_2}$).

\vspace{0.5cm}

\begin{example}\label{caps}
The SFT-potential of $\mathbb{C}/\mathbb{Z}_2$ has the form
$$\mathbf{F}_{2}(t_0,t_1,p)=\frac{1}{\hbar}\left[At_0t_1^2+Bt_1^4+Ct_1p_1+Dp_2\right]$$
and a comparison with the genus $0$ GW-potential of $\mathbb{P}^1_{2,2}$ (see for instance \cite{DZ1})
$$\textstyle{\mathbf{F}_{2,2}=\frac{1}{\hbar}\left[\frac{1}{2}t_0^2s+\frac{1}{4}(t_{(1,1)}^2+t_{(1,2)}^2)t_0- \frac{1}{96}(t_{(1,1)}^4+t_{(1,2)}^4)+t_{(1,1)}t_{(1,2)}\mathrm{e}^{s}z+\frac{1}{2}\mathrm{e}^{2s}z^2\right]}$$
(where the second index of the cohomology variables refers to one of the two orbifold points) promptly gives $$A=\frac{1}{4}\qquad B=-\frac{1}{96}\qquad C=1\qquad D=\frac{1}{2}$$

\noindent Similarly one computes
$$\textstyle{\mathbf{F}_{3}=\frac{1}{\hbar}\left[\frac{1}{3}t_0t_1t_2+\frac{1}{18}t_1^3-\frac{1}{36}t_1^2t_2^2 +\frac{1}{648}t_1t_2^4-\frac{1}{19440}t_3^6+\left(t_1+\frac{1}{6}t_2^2\right)p_1+ \frac{1}{2}t_2p_2+\frac{1}{3}p_3\right]}$$
\vspace{0.3cm}
\begin{equation*}
\begin{split}
&\textstyle{\mathbf{F}_{4}=\frac{1}{\hbar}\left[-\frac{t_3^8}{4128768}+\frac{t_2 t_3^6}{73728}-\frac{t_1 t_3^5}{30720}-\frac{t_2^2 t_3^4}{3072}+\frac{1}{384} t_1 t_2
   t_3^3+\frac{1}{384} t_2^3 t_3^2-\frac{1}{64} t_1^2 t_3^2-\frac{1}{32} t_1 t_2^2 t_3+\frac{1}{4} t_0 t_1
   t_3\right.}\\
&\textstyle{\hspace{1cm}\left.-\frac{t_2^4}{192}+\frac{1}{8}
   t_0 t_2^2+\frac{1}{8} t_1^2 t_2+
   \left(\frac{t_3^3}{96}+\frac{1}{4} t_2 t_3+t_1\right)p_1 +
   \left(\frac{t_3^2}{8}+\frac{t_2}{2}\right)p_2+\frac{1}{3} t_3 p_3+\frac{1}{4} p_4\right]}
\end{split}
\end{equation*}
\vspace{0.3cm}
\begin{equation*}
\begin{split}
&\textstyle{\mathbf{F}_{5}=\frac{1}{\hbar}\left[-\frac{7
   t_4^{10}}{8100000000}+\frac{7 t_3 t_4^8}{90000000}-\frac{t_2 t_4^7}{3150000}-\frac{13 t_3^2 t_4^6}{4500000}+\frac{t_1
   t_4^6}{2250000}+\frac{t_3^5}{3000}+\frac{11 t_2 t_3 t_4^5}{375000}+\frac{7 t_3^3 t_4^4}{150000}\right.}\\
&\textstyle{\hspace{1cm}-\frac{t_2^2
   t_4^4}{7500}-\frac{t_1 t_3 t_4^4}{15000}-\frac{1}{150} t_1 t_3^3-\frac{t_2 t_3^2 t_4^3}{1500}+\frac{1}{750} t_1 t_2
   t_4^3+\frac{1}{10} t_1 t_2^2-\frac{1}{50} t_2^2 t_3^2-\frac{3 t_3^4 t_4^2}{10000}-\frac{1}{100} t_1^2 t_4^2}\\
&\textstyle{\hspace{1cm}+\frac{1}{500} t_1
   t_3^2 t_4^2+\frac{1}{250} t_2^2 t_3 t_4^2+\frac{1}{10} t_1^2 t_3+\frac{1}{5} t_0 t_2 t_3-\frac{1}{150} t_2^3 t_4+\frac{1}{250}
   t_2 t_3^3 t_4+\frac{1}{5} t_0 t_1 t_4-\frac{1}{25} t_1 t_2 t_3 t_4}\\
&\textstyle{\hspace{1cm}+ \left.\left(\frac{t_4^4}{3000}+\frac{1}{50} t_3 t_4^2+\frac{t_2 t_4}{5}+\frac{t_3^2}{10}+t_1\right)p_1+
   \left(\frac{t_4^3}{75}+\frac{t_3 t_4}{5}+\frac{t_2}{2}\right)p_2+
   \left(\frac{t_4^2}{10}+\frac{t_3}{3}\right)p_3 +\frac{1}{4} t_4 p_4+\frac{1}{5}p_5\right]
}
\end{split}
\end{equation*}
\end{example}

Now we can proceed with our gluing procedure, which consists, according to Theorem $2.5.3$ of \cite{EGH}, in obtaining the GW-potential for the genus $g'$ orbicurve $S_{\alpha_1,\ldots,\alpha_a}=\Sigma_{g',(z_1,\alpha_1),\ldots,(z_a,\alpha_a)}$, with orbifold points $z_1,\ldots,z_a$ and local groups $\mathbb{Z}_{\alpha_1},\ldots,\mathbb{Z}_{\alpha_a}$ as
\begin{equation}
\mathbf{F}_{S_{\alpha_1,\ldots,\alpha_a}}(t,s_2)= \mathrm{log}\left[\mathrm{exp}(\overrightarrow{\mathbf{F}}_{\alpha_1}(t,p))\ldots \mathrm{exp}(\overrightarrow{\mathbf{F}}_{\alpha_a}(t,p)) \, \mathrm{exp}(\mathbf{F}_S(t_0,s_2,q))\right]
\end{equation}
which gives, according for our general formula for the potential of the orbifold cap, the following expression of the genus $g$ orbifold Gromov-Witten potential of the orbicurve $S_{\alpha_1,\ldots,\alpha_a}$
\begin{equation}\label{potential}
\begin{split}
\mathbf{F}_{S_{\alpha_1,\ldots,\alpha_a}}=&\frac{1}{\hbar}\frac{t_0^2s_2}{2}-\frac{s_2}{24}\\
&+\sum_{r=1}^a \sum_{\substack{i_{(0,r)},\ldots,i_{(\alpha_r-1,r)} \\ \sum(\alpha_r-k)i_{(k,r)}=2\alpha_r}} A^r_{i_{(0,r)}\ldots i_{(\alpha_r-1,r)}}t_0^{(i_0,r)} t_{(1,r)}^{i_{(1,r)}}\ldots t_{(\alpha_r-1,r)}^{i_{(\alpha_r-1,r)}}\\
&+\sum_{g=0}^{\infty}\sum_{d=0}^{\infty}\sum_{|\mu_1|,\ldots,|\mu_a|=d}\mathrm{e}^{d s_2}H^{\Sigma_{g'}}_{g,d}(\mu_1,\ldots,\mu_a)\\
&\prod_{r=1}^{a}\left(\sum_{i_{(1,r)},\ldots,i_{(\alpha_r-1,r)} }B^r_{i_{(1,r)}\ldots i_{(\alpha_r-1,r)}}t_{(1,r)}^{i_{(1,r)}}\ldots t_{(\alpha_r-1,r)}^{i_{(\alpha_r-1,r)}}\right)^{\mu_r}\hbar^{g-1}
\end{split}
\end{equation}
where, as in Example \ref{caps}, the second index in round brackets, as well as the upper index of $A$ and $B$,  specifies one of the $a$ orbifold points.

Notice here that, in the above formula, the sum over the branching configurations $\mu_1,\ldots,\mu_a$, where $\mu_k=(\mu_{1,k},\ldots,\mu_{s_k,k})$, $k=1,\ldots,a$ is the local branching degree over the $k$-th point, involves only those terms for which each of the $\mu_{j,k}$ is less or equal than $\alpha_k$. This phenomenon is very important and sometimes it leaves only a finite number of nonzero terms in genus $0$, making the rational GW-potential $\mathbf{f}_{S_{\alpha_1,\ldots,\alpha_a}}:=\left.\left(\hbar\, \mathbf{F}_{S_{\alpha_1,\ldots,\alpha_a}}\right)\right|_{\hbar=0}$ a polynomial in the variables $t_0,t_{(1,1)},\ldots,t_{(\alpha_1-1,1)},\ldots,t_{(1,a)},\ldots,t_{(\alpha_a-1,a)},\mathrm{e}^{s_2}$. Our next task will consist in classifying exactly these cases, at least for $g'=0$.

\section{Polynomial $\mathbb{P}^1$-orbifolds}
 
In this section we study the rational orbifold Gromov-Witten theory of those genus $0$ orbicurves which give rise to a polynomial quantum cohomology, i.e. those whose associated genus $0$ Gromov-Witten potential is polynomial in the variables $t_0,t_{(1,1)},\ldots,t_{(\alpha_1-1,1)},\ldots,$ $t_{(1,a)}, \ldots,t_{(\alpha_a-1,a)},\mathrm{e}^s$. We will call such orbifold Riemann surfaces \emph{polynomial $\mathbb{P}^1$-orbifolds}.

\begin{lemma}\label{polyP1}
The only polynomial $\mathbb{P}^1$-orbifolds are those $\mathbb{P}^1_{(z_1,\alpha_1),\ldots,(z_a,\alpha_a)}$ such that:
\begin{itemize}
\item[1)] $a=0,1,2$, for any $(\alpha_1,\alpha_2)$: $$\mathbb{P}^1_{\alpha_1,\alpha_2}$$
\item[2)] $a=3$, $(\alpha_1,\alpha_2,\alpha_3)=(2,2,l-2)$, $l\ge 4$: $$\mathbb{P}^1_{2,2,l-2}$$
\item[3)] $a=3$, $(\alpha_1,\alpha_2,\alpha_3)=(2,3,3),(2,3,4),(2,3,5)$:
$$\mathbb{P}^1_{2,3,5} \qquad \mathbb{P}^1_{2,3,4} \qquad \mathbb{P}^1_{2,3,5}$$
\end{itemize}
\end{lemma}
\begin{proof}
A necessary and sufficient condition for the finiteness of the potential is obtained by using Riemann-Hurwitz theorem. Indeed, validity of the Riemann-Hurwitz relation
$$\sum_{k=1}^a \sum_{i=1}^{s_k}(\mu_{i,k}-1)=2d-2$$
is a necessary condition for existence of a degree $d$ covering with branching profile given by $\mu_1,\ldots,\mu_a$, with $\mu_k=(\mu_{1,k},\ldots,\mu_{s_k,k})$. As we already noticed, the branched coverings appearing in our formula for the potential of $\mathbb{P}^1_{(z_1,\alpha_1),\ldots,(z_a,\alpha_a)}$ have a branching profile with local branching degrees $\mu_k=(\mu_{1,k},\ldots,\mu_{s_k,k})$, $k=1,\ldots,a$ such that $\mu_{i,k}\le\alpha_k$, $i=1,\ldots,s_k$. Denoting by $n_{j,k}$ the number of occurrences of $j$ in the set $\{\mu_{1,k},\ldots,\mu_{s_k,k}\}$, we can then rewrite the Riemann-Hurwitz relation as
$$\sum_{k=1}^a \sum_{j=1}^{\alpha_k} (j-1)n_{j,k}=2d-2$$
We are interested in singling out the situations where there are only a finite number of positive integer solutions $n_{j,k}$ to this equation, together with the constant degree conditions $d=\sum_{j=1}^{\alpha_k}j\,n_{j,k}$, $k=1.\ldots,a$. It is immediately clear (consider the branching configurations where $n_{1,1}=4,$ $n_{2,1}=2r-2$, $n_{2,2}=n_{2,3}=n_{2,4}=2r$ are the only nonzero coefficients) that this never happens when $a\ge 4$, while it is always the case when $a\le 2$. Finally, when $a=3$ we distinguish between the following cases. For $\alpha_1=\alpha_2=2$, to ensure positiveness of $n_{1,1},n_{2,1},n_{1,2}$ we must have
$$(\alpha_3+2)n_{\alpha_3,3}+(\alpha_3+1)n_{\alpha_3-1,3}+\ldots+3n_{1,3}-4\le2n_{2,2}$$
$$(\alpha_3+1)n_{\alpha_3,3}+\alpha_3n_{\alpha_3-1,3}+\ldots+2n_{1,3}-2\ge n_{2,2}$$
$$\alpha_3n_{\alpha_3,3}+(\alpha_3-1)n_{\alpha_3-1,3}+\ldots+n_{1,3}\ge2n_{2,2}$$
whence $2(n_{\alpha_3,3}+\ldots+n_{1,3})\le 4$, so we get a finite number of integer solutions for any $\alpha_3$. For $\alpha_1=2$ and $\alpha_2=3$ we have similarly
$$(\alpha_3+2)n_{\alpha_3,3}+(\alpha_3+1)n_{\alpha_3-1,3}+\ldots+3n_{1,3}-4\le2n_{2,2}+n_{3,1}$$
$$(\alpha_3+1)n_{\alpha_3,3}+\alpha_3n_{\alpha_3-1,3}+\ldots+2n_{1,3}-2\ge n_{2,2}+2n_{3,1}$$
$$\alpha_3n_{\alpha_3,3}+(\alpha_3-1)n_{\alpha_3-1,3}+\ldots+n_{1,3}\ge2n_{2,2}$$
whence $(6-\alpha_3)n_{\alpha_3,3}+\ldots+5n_{1,3}\le 12$, so the solutions are in a finite number if $\alpha_3\le 5$. Vice versa, if $\alpha_3\ge6$, we have infinite sequences of solutions like $n_{2,2}=3r$, $n_{3,1}=2r-2$, $n_{6,3}=r$ and $n_{j,k}=0$ otherwise.\\
Other values of $\alpha_1, \alpha_2$ always give infinite solutions.

\vspace{0.3cm}

Notice that the Riemann-Hurwitz relation is not a sufficient condition and the problem of determining all the admissible (satisfying Riemann-Hurwitz) covering configurations that are not actually geometrically realizable is classical and still open (see e.g. \cite{PP} for a review about this challenging topic). Nonetheless the non-realizable cases are quite exceptional and (as shown in \cite{PP}) by no means influence infiniteness of the number of actual coverings. This ensures that we have actually found all the polynomial $\mathbb{P}^1$-orbifolds.
\end{proof}

In the above classification, case $1)$ was studied in \cite{MT}, where the quantum cohomology of $\mathbb{P}^1_{\alpha_1,\alpha_2}$ was shown to be isomorphic, as Frobenius manifold, to the space $M(A_l,k)$ of Fourier polynomials invariant with respect to the extended affine Weyl group of the root system $A_l$, with the choice of the $k$-th root and with $l=\alpha_1+\alpha_2-1$ and $k=\alpha_1$. We will not remind here the general procedure (see \cite{DZ1}) to construct the Frobenius structure on the space $M(R,k)$. However we recall that polynomial Frobenius manifolds of dimension $l+1$ can be associated with root systems of type $D_l$ and $E_l$, but in this case the choice of the $k$-th root is forced to be at the bifurcation of the Dynkin diagram. Moreover, the relevant Frobenius potential happens to be homogeneous with respect to a specific grading for the coordinates, and this grading coincides (up to an irrelevant factor) with the one we defined above once we compare the manifolds $M(D_l,l-2)$ with $QH^*_\mathrm{orb}(\mathbb{P}^1_{2,2,l-2})$ and $M(E_l,4)$ with $QH^*_\mathrm{orb}(\mathbb{P}^1_{2,3,l-3})$.

\vspace{0.5cm}

With the idea to investigate this correspondence, we plan to compute explicitly the genus $0$ Gromov-Witten potential $\mathbf{f}_{\mathbb{P}^1_{\alpha_1,\alpha_2,\alpha_3}}$ of the above polynomial $\mathbb{P}^1$-orbifolds. This can be done easily (although computations can get quite cumbersome) by using equation (\ref{potential}) above, together with Riemann-Hurwitz relation, as in the proof of the above Lemma, to predict all the (a priori) non-zero Hurwitz numbers appearing as coefficients in our formula. This leaves us with a closed expression for $\mathbf{f}_{\mathbb{P}^1_{\alpha_1,\alpha_2,\alpha_3}}$ which still contains a finite number of Hurwitz coefficients. These can then be determined by imposing WDVV equations. Without reporting the details of these computations, we give some examples of the results of this construction.

\begin{example}\label{examples}
\begin{equation*}
\begin{split}
&\textstyle{\mathbf{f}_{\mathbb{P}^1_{2,2,2}}= \frac{1}{4} e^{4 t_4} z^4+\frac{1}{2} e^{2 t_4} \left(t_1^2+t_2^2+t_3^2\right) z^2+e^{t_4} t_1 t_2 t_3 z+\frac{1}{4} t_0
   \left(t_1^2+t_2^2+t_3^2\right)+\frac{1}{96} \left(t_1^4+t_2^4+t_3^4\right)}\\
&\textstyle{\hspace{1.4cm}   +\frac{1}{2} t_0^2 t_4}
\end{split}
\end{equation*}

\begin{equation*}
\begin{split}
&\textstyle{\mathbf{f}_{\mathbb{P}^1_{2,2,3}}= \frac{1}{6} e^{6 t_5} z^6+\frac{1}{4} e^{4 t_5} t_4^2 z^4+e^{3 t_5} t_1 t_2 z^3+\frac{1}{2} e^{2 t_5}
   \left(\frac{t_4^2}{6}+t_3\right)^2 z^2+\frac{1}{2} e^{2 t_5} \left(t_1^2+t_2^2\right) t_4 z^2}\\
&\textstyle{\hspace{1.4cm} +e^{t_5} t_1 t_2
\left(\frac{t_4^2}{6}+t_3\right) z-\frac{t_4^6}{19440}+\frac{1}{648} t_3 t_4^4+\frac{t_3^3}{18}-\frac{1}{36} t_3^2
   t_4^2+\frac{1}{4} t_0 \left(t_1^2+t_2^2\right)-\frac{1}{96} \left(t_1^4+t_2^4\right)}\\
&\textstyle{\hspace{1.4cm}   +\frac{1}{3} t_0 t_3 t_4+\frac{1}{2}
   t_0^2 t_5}
\end{split}
\end{equation*}

\begin{equation*}
\begin{split}
&\textstyle{\mathbf{f}_{\mathbb{P}^1_{2,2,4}}= \frac{1}{8} e^{8 t_6} z^8+\frac{1}{6} e^{6 t_6} t_5^2 z^6+e^{4 t_6} \left(\frac{t_5^2}{8}+\frac{t_4}{2}\right)^2 z^4+\frac{1}{2}
   e^{4 t_6} \left(t_1^2+t_2^2\right) z^4+e^{3 t_6} t_1 t_2 t_5 z^3}\\
&\textstyle{\hspace{1.4cm}    +\frac{1}{2} e^{2 t_6} \left(\frac{t_5^3}{96}+\frac{t_4
   t_5}{4}+t_3\right)^2 z^2+e^{2 t_6} \left(t_1^2+t_2^2\right) \left(\frac{t_5^2}{8}+\frac{t_4}{2}\right) z^2+e^{t_6} t_1 t_2
   \left(\frac{t_5^3}{96}+\frac{t_4 t_5}{4}+t_3\right) z}\\
&\textstyle{\hspace{1.4cm}     -\frac{t_5^8}{4128768}+\frac{t_4 t_5^6}{73728}-\frac{t_3
   t_5^5}{30720}-\frac{t_4^4}{192}-\frac{t_4^2 t_5^4}{3072}+\frac{1}{384} t_3 t_4 t_5^3+\frac{1}{8} t_0 t_4^2+\frac{1}{384}
   t_4^3 t_5^2-\frac{1}{64} t_3^2 t_5^2}\\
&\textstyle{\hspace{1.4cm}       +\frac{1}{4} t_0 \left(t_1^2+t_2^2\right)+\frac{1}{96}
   \left(-t_1^4-t_2^4\right)+\frac{1}{8} t_3^2 t_4-\frac{1}{32} t_3 t_4^2 t_5+\frac{1}{4} t_0 t_3 t_5+\frac{1}{2} t_0^2 t_6}\\
\end{split}
\end{equation*}

\begin{equation*}
\begin{split}
&\textstyle{\mathbf{f}_{\mathbb{P}^1_{2,3,3}}=-\frac{t_5^4}{96}+\frac{1}{3} e^{3 t_6} t_5^3+\frac{1}{2} e^{6 t_6} t_5^2+\frac{1}{4} t_0 t_5^2+\frac{1}{2} e^{2 t_6} t_2 t_4 t_5^2+e^{5 t_6} t_2 t_4 t_5+e^{t_6}
   \left(\frac{t_2^2}{6}+t_1\right) \left(\frac{t_4^2}{6}+t_3\right) t_5}\\
&\textstyle{\hspace{1.4cm}+e^{3 t_6} \left(t_2 \left(\frac{t_2^2}{6}+t_1\right)+t_4
   \left(\frac{t_4^2}{6}+t_3\right)\right) t_5+\frac{e^{12 t_6}}{12}+\frac{1}{4} e^{4 t_6} t_2^2 t_4^2+\frac{1}{18} \left(t_1^3+t_3^3\right)+\frac{1}{2} e^{8 t_6} t_2
   t_4}\\
&\textstyle{\hspace{1.4cm}+\frac{1}{3} t_0 \left(t_1 t_2+t_3 t_4\right)+e^{4 t_6} \left(\frac{t_2^2}{6}+t_1\right) \left(\frac{t_4^2}{6}+t_3\right)+\frac{1}{36} \left(-t_1^2 t_2^2-t_3^2
   t_4^2\right)+\frac{1}{6} e^{6 t_6} \left(t_2^3+t_4^3\right)}\\
&\textstyle{\hspace{1.4cm}+\frac{1}{648} \left(t_1 t_2^4+t_3 t_4^4\right)+\frac{-t_2^6-t_4^6}{19440}+\frac{1}{2} e^{2 t_6}
   \left(t_4 \left(\frac{t_2^2}{6}+t_1\right)^2+t_2 \left(\frac{t_4^2}{6}+t_3\right)^2\right)+\frac{1}{2} t_0^2 t_6}
\end{split}
\end{equation*}
\end{example}

All the above potentials turn out to coincide with know examples from extended affine Weyl groups and, in particular, one can directly verify the isomorphism of Frobenius manifolds
$$M(E_l,4)\cong QH^*_\mathrm{orb}(\mathbb{P}^1_{2,3,l-3})\qquad l=6,7,8$$
together with instances of the other isomorphisms appearing in theorem \ref{mirror}. However a complete proof of our mirror result will require a less computational approach, which is what we plan for the next sections.

\section{Space of tri-polynomials}

We will denote by $M_{p,q,r}$ the space of polynomials (we will refer to them as tri-polynomials) of the form
$$F(x,y,z)=-xyz+P_1(x)+P_2(y)+P_3(z)$$
where
$$P_1(x)=\sum_{k=1}^{p}a_k x^k \hspace{2cm} P_2(y)=\sum_{k=1}^{q}b_k y^k \hspace{2cm} P_3(z)=\sum_{k=0}^{r}c_k (\mathrm{e}^{d} z)^k$$
and normalized by $a_p=b_q=c_r=1$.  This means that $M_{p,q,r}\cong\mathbb{C}^{p+q+r-2}\times\mathbb{C}^*$.
It is an easy exercise with generators and relations to show that the condition on $p,q,r$ such that the local algebra $\mathbb{C}[x,y,z]/J_F$, where $J_F$ is the Jacobian ideal of $F$ (i.e. $J_F=(\partial_x F,\partial_y F, \partial_z F)$), is isomorphic as a $\C$-module to the tangent space $T_FM_{p,q,r}$ is $\frac{1}{p}+\frac{1}{q}+\frac{1}{r}>1$. We will assume this condition is always verified in what follows. Moreover, up tp permutation of $x,y,z$, we will just need to consider the cases
\begin{itemize}
 \item[($A$)] $(p,q,r)=(p,q,1),\hspace{1cm}p,q=1,2,\ldots$
 \item[($D$)] $(p,q,r)=(2,2,r),\hspace{1cm}r=2,3,\ldots$
 \item[($E$)]$(p,q,r)=(2,3,r),\hspace{1cm}r=3,4,5$
\end{itemize}

We are going to define a Frobenius manifold structure on $M_{p,q,r}$. Similarly to the case of Laurent polynomials (see e.g. \cite{D},\cite{MT}), the isomorphism of each tangent space $T_{F}M_{p,q,r}$ with the local algebra $\mathbb{C}[x,y,z]/J_F$, defines a commutative associative algebra structure on $T_{F}M_{p,q,r}$, with unity $e=\partial / \partial c_0$.

We then equip each tangent space with the residue pairing
\begin{equation}\label{metric}
(\partial,\partial ')_F=\underset{\partial_x F=\partial_y F=\partial_z F=0}{\underset{F(x,y,z)\neq\infty}{\mathrm{res}}} \frac{\partial(F)\, \partial'(F)}{\partial_x F\, \partial_y F\, \partial_z F}\, \di x \wedge \di y \wedge \di z
\end{equation}
and assign the following grading to the variables
$$\mathrm{deg}\,a_i=-2+\frac{2i}{p} \hspace{1cm} i=1,\ldots,p-1$$
$$\mathrm{deg}\,b_j=-2+\frac{2j}{q} \hspace{1cm} j=1,\ldots,q-1$$
$$\mathrm{deg}\,c_k=-2+\frac{2k}{r} \hspace{1cm} k=0,\ldots,r-1$$
$$\mathrm{deg}\,d = 0$$
which gives the following Euler vector field
$$E=\sum_{i=1}^{p-1} (1-\frac{i}{p})a_i\frac{\partial}{\partial a_i}+\sum_{j=1}^{q-1} (1-\frac{j}{q})b_j\frac{\partial}{\partial b_j} + \sum_{k=0}^{r-1} (1-\frac{k}{r})c_k\frac{\partial}{\partial c_k}+(-1+\frac{1}{p}+\frac{1}{q}+\frac{1}{r})\frac{\partial}{\partial d}$$
Notice that, as usual in singularity theory, $E$ corresponds to $F(x,y,z)$ via the isomorphism $TM_{p,q,r}\simeq \C[x,y,z]/J_F$.

All these structures are easily verified to be compatible, hence it only remains to prove flatness of the metric (\ref{metric}) and potentiality (i.e. the symmetry of $\nabla_{X} g(Y,Z\circ W)$ in $X,Y,Z,W$). The first requirement is achieved by exhibiting a system of flat coordinates. We start with the cases $A$ and $D$. Let
$$\alpha_i=\underset{x=\infty}{\mathrm{res}} \frac{p}{i} \frac{F(x,0,0)^{1-\frac{i}{p}}}{x}\di x \hspace{1cm} i=1,\ldots,p-1$$
$$\beta_j=\underset{y=\infty}{\mathrm{res}} \frac{q}{j} \frac{F(0,y,0)^{1-\frac{j}{q}}}{y}\di y \hspace{1cm} j=1,\ldots,q-1$$
$$\gamma_k=\underset{z=\infty}{\mathrm{res}} \frac{r}{k} \frac{F(0,0,z)^{1-\frac{k}{r}}}{\sqrt{z^2-4}}\di z \hspace{1cm} k=0,\ldots,r-1$$
where, as specified, either $(p,q,r)=(p,q,1)$ or $(p,q,r)=(2,2,r)$.

\begin{lemma}\label{flat}
In the cases $A$ and $D$, the functions $\alpha_i,\beta_j,\gamma_k:\,M_{p,q,r}\to\C$, together with $d:\, M_{p,q,r}\to\C$, form a system of flat coordinates for the metric (\ref{metric}), such that the only nonzero pairings of basis vectors are given by
$$(\partial_{\alpha_{i_1}},\partial_{\alpha_{i_2}})=\frac{1}{p}\,\delta_{i_1+i_2,p}$$
$$(\partial_{\beta_{j_1}},\partial_{\beta_{j_2}})=\frac{1}{q}\,\delta_{j_1+j_2,q}$$
$$(\partial_{\gamma_{k_1}},\partial_{\gamma_{k_2}})=\frac{1}{r}\,\delta_{k_1+k_2,r}$$
$$(\partial_{\gamma_0},\partial_{d})=1$$
with indices ranging as above.
\end{lemma}
\begin{proof}
In the $A$ case the formulas for $\alpha_i$ and $\beta_j$ coincide with the ones given in \cite{DZ1} for the space of Laurent polynomials, hence only the $D$ case is left to prove.

From the above formulae we see that $\alpha_1=a_1$ and $\beta_1=b_1$. We now prove that $(\partial_{\gamma_{k_1}},\partial_{\gamma_{k_2}})=\frac{1}{r}\,\delta_{k_1+k_2,r}$. Consider the solutions, with respect to $z$, to $F(x,y,z)=\lambda$, for $\lambda$ and $z$ near infinity, and denote it by $z=z(\lambda;x,y)$. Using the chain rule we get $\partial_{\gamma_k}F=-(\partial_{z}F)(\partial_{\gamma_k}z)$. Moreover, in the $D$ case, the residues for the matrix elements of the metric with respect to $x$ and $y$ localize at the poles
\begin{equation*}
\left\{ \begin{array}{l} \partial_x F(x,y,z)=0\\ \partial_y F(x,y,z)=0 \end{array} \right. \hspace{1cm} \Rightarrow \hspace{1cm}
\left\{ \begin{array}{l} x_p=\frac{b_1 z+2a_1}{z^2-4}\\ y_p=\frac{a_1z+2b_1}{z^2-4} \end{array}\right.
\end{equation*}
Hence, applying the residue formula in many variables gives
\begin{equation*}
\begin{split}
(\partial_{\gamma_{k_1}},\partial_{\gamma_{k_2}})&=-\underset{\lambda=\infty}{\mathrm{res}} \frac{(\partial_{\gamma_{k_1}}z)(\partial_{\gamma_{k_2}}z)}{z^2-4}\,\di\lambda\\
&=-\underset{\lambda=\infty}{\mathrm{res}} (\partial_{\gamma_{k_1}}\mathrm{log}(z+\sqrt{z^2-4})) (\partial_{\gamma_{k_2}}\mathrm{log}(z+\sqrt{z^2-4}))\,\di\lambda
\end{split}
\end{equation*}
and the coordinates $(x_p,y_p)$ make no appearance in the formula, thanks to the fact that the coordinates $\gamma_i$ only depend on $c_i$ and $d$.

The result then follows if we consider the expansion, for $\lambda$ near infinity, of $\mathrm{log}(z+\sqrt{z^2-4})$) where the coefficient of $\lambda^{k/r}$ is given by $\gamma_k$. Indeed, let the coefficients $\gamma_k$ be defined by
$$\mathrm{log}(z+\sqrt{z^2-4})=\frac{1}{r}\left[\mathrm{log}\lambda- \gamma_{r-1}\lambda^{-\frac{1}{r}}-\ldots-\gamma_1\lambda^{-\frac{r-1}{r}}-\gamma_0\lambda^{-1}\right]+O(\lambda^{-1-\frac{1}{r}})$$
Then one has
\begin{equation*}
 \underset{z=\infty}{\mathrm{res}} \frac{r}{k}\, \frac{\lambda^{1-\frac{k}{r}}}{\sqrt{z^2-4}}\di z =
 \underset{z=\infty}{\mathrm{res}} \frac{r}{k}\,  \lambda^{1-\frac{k}{r}} \left(\frac{\partial}{\partial \lambda^{\frac{1}{r}}}\left(\mathrm{log}\frac{z+\sqrt{z^2-4}}{\lambda^{\frac{1}{r}}}\right) +\frac{1}{\lambda^{\frac{1}{r}}}\right)\di \lambda^{\frac{1}{r}}=\gamma_k
\end{equation*}

Expressions for $(\partial_{\alpha_{1}},\partial_{\alpha_{1}})$, $(\partial_{\beta_{1}},\partial_{\beta_{1}})$ and $(\partial_{\gamma_0},\partial_{d})$ are proved analogously. This also shows that these are the only nonzero entries.
\end{proof}

The $E$ case is to be dealt with similarly, but it needs more care. First one finds the poles with respect to $x$ and $y$ in the integrand of (\ref{metric}):
\begin{equation*}
\left\{ \begin{array}{l} x_{p}=\frac{1}{24}\left(z^3-12a_1-4b_2z\pm z\sqrt{z^4-8b_2z^2-24a_1z+16b_2^2-48b_1}\right)\\ y_p=\frac{1}{12}\left(z^2-4b_2\pm\sqrt{z^4-8b_2z^2-24a_1z+16b_2^2-48b_1}\right) \end{array}\right.
\end{equation*}
(we will use subscripts $p=p_1$ for the choice of plus and $p=p_2$ for the choice of minus in the coordinates above) and, applying once more the residue formula, for any two coordinates $t_1$ and $t_2$ one gets
\begin{equation*}
(\partial_{t_1},\partial_{t_2})=\underset{\lambda=\infty}{\mathrm{res}} \left( \frac{\left[(\partial_{t_1}z)(\partial_{t_2}z)\right]_{(x=x_{p_1},y=y_{p_1})}- \left[(\partial_{t_1}z)(\partial_{t_2}z)\right]_{(x=x_{p_2},y=y_{p_2})}}{\sqrt{z^4-8b_2z^2-24a_1z+16b_2^2-48b_1}}\right)\,\di\lambda
\end{equation*}
At this point the esasiest thing is probably just computing the flat coordinates from the topological side (recall that we are able to find the explicit expression of the Frobenius potential for $\mathbb{P}^1_{2,3,r}$, $r=3,4,5$) and plug them into the above formula to check that they are flat coordinates for the metric (\ref{metric}) too.
For instance, for the $E_6$ case (see Example \ref{examples}), one gets:
$$\alpha_1=a_1-8\mathrm{e}^{3d}$$
$$\beta_1=b_1-\frac{b_2^2}{6}+3c_2\mathrm{e}^{2d}\ ,\hspace{1cm} \beta_2=b_2$$
$$\gamma_0=c_0+2a_2b_2\mathrm{e}^{2d}+6a_1\mathrm{e}^{3d}-18\mathrm{e}^{6d}$$
$$\gamma_1=c_1-\frac{c_2^2}{6}+3b_2\mathrm{e}^{2d}\ ,\hspace{1cm} \gamma_2=b_2$$

As for potentiality, one can use the same technique as in Hertling's book on singularities and Frobenius structures \cite{H}. Recall that an $F$-manifold structure $(M,\circ,e)$ on a complex manifold $M$ is given by a commutative and associative multiplication $\circ$ on $TM$ and a unity vector field $e$ such that $\mathrm{Lie}_{X\circ Y}(\circ) = X\circ \mathrm{Lie}_{Y}(\circ)+ Y\circ \mathrm{Lie}_{X}(\circ)$. Indeed we use the following result from \cite{H}.

\begin{thm}[\cite{H}]
Let $(M,\circ,e,g)$ be a manifold with a commutative and associative multiplication $\circ$ on $TM$, a unity vector field $e$ and a metric $g$ which is multiplication invariant. Denote by $\epsilon$ the $1$-form (coidentity) $g(\cdot,e)$. Then the following are equivalent:
\begin{itemize}
\item[i)] $M$ carries a structure of an $F$-manifold and $\epsilon$ is closed
\item[ii)] $\nabla_{\cdot}(\cdot,\cdot\circ\cdot)$ is symmetric in all four arguments
\end{itemize}
\end{thm}

In our case $e=\partial_{c_0}=\partial_{\gamma_0}$, so $\di\epsilon=0$ follows from flatness of $e$. The $F$-manifold condition can be deduced from the following theorem.

\begin{thm}[\cite{H}]
Let $\pi_{Z}:Z \twoheadrightarrow M$ be a submersion between the manifolds $Z$ and $M$, with $\mathrm{dim}\,M=n$. Let $C\subset Z$ be a an $n$-dimensional reduced subvariety such that the restriction $\pi_{C}:C \twoheadrightarrow M$ is finite. Finally, let $\alpha_Z$ be a $1$-form on $Z$ such that, for any local lift $\tilde{X}\in \Gamma(TZ)$ of the zero vector field on $M$, $\alpha_Z(\tilde{X})|_C=0$.

Then the map $$\Gamma(TM)\to(\pi_C)_*\mathcal{O}_C \hspace{1cm} X\mapsto \alpha_Z(\tilde{X})|_C$$ is well defined, and provides $M$ with a structure of $F$-manifold (with generically semi-simple multiplication) if and only if it is an isomorphism and $\alpha_Z|_{C_{\mathrm{reg}}}$ is exact.
\end{thm}

In our case $Z=\mathbb{C}^3\times M_{p,q,r}$, $M=M_{p,q,r}$, $C=\{\di_{\mathbb{C}^3}F(x,y,z)=0\}$ so that $(\pi_C)_*\mathcal{O}_C=(\C[x,y,z]\otimes \mathcal{O}_M)/J_F$, and $\alpha_Z=\di_Z F$, so the hypothesis of the above theorem are easily verified.

Notice that the spaces of tri-polynomials given by $M_{p,q,1}$ are trivially isomorphic, as Frobenius manifolds, to the spaces of Laurent polynomials $M_{p,q}$ considered by Dubrovin and Zhang (\cite{DZ1}) and Milanov and Tseng (\cite{MT}).

\section{The mirror theorem}

In this section we prove the isomorphism between $M_{p,q,r}$ and $QH_{\mathrm{orb}}^*(\mathbb{P}^1_{p,q,r})$, for $\frac{1}{p}+\frac{1}{q}+\frac{1}{r}>1$, as Frobenius manifolds. Actually, after the work of Milanov and Tseng (\cite{MT}), only the cases $M_{2,2,r}$ and the exceptional cases $M_{2,3,3}$, $M_{2,3,4}$, $M_{2,3,5}$ are left to be proven. In order to do that, we are going to show that we can find a point $m\in M_{p,q,r}$ (and hence a dense subset) of our Frobenius manifolds where the linear operator $U: T_mM_{p,q,r}\to T_mM_{p,q,r}$ of multiplication by $E_m$ has pairwise distinct eigenvalues. This ensures (see e.g. \cite{D1}) that the Frobenius manifold is semisimple (meaning that the Frobenius algebra is semisimple on a dense subset). Moreover, if we prove that the algebras $T_mM_{p,q,r}$ and $T_mQH^*_{\mathrm{orb}}(\mathbb{P}^1_{p,q,r})$ at $m$ are isomorphic, then the isomorphism of Frobenius structures follows from the identification of $e$, $E$ and the metric, via the diffeomorphism induced by sending the flat coordinates on $M_{p,q,r}$ to the components of homogeneous basis of cohomology classes in $QH^*_\mathrm{orb}(\mathbb{P}^1_{p,q,r})$. In fact (see \cite{D1}, Lemma $3.3$), knowing metric, unity vector field $e$ and Euler vector field $E$, together with the algebra structure at a point $m$ where $U=(E\bullet_m)$ has pairwise distinct eigenvalues, allows one to reconstruct uniquely the whole Frobenius structure.

After having established the algebra isomorphism, we are going to compute the operator $U$ in the linear basis given by flat coordinates, at the origin $0$ of the flat coordinate system. Here eigenvalues can be calculated explicitly but they fail to be distinct. However, in the $M_{2,2,r}$ case, we are able to write down an explicit expression for $U$ along a two-dimensional submanifold of $M_{p,q,r}$ and, by using perturbation theory along this submanifold, we show that degeneracy of the eigenvalues is completely removed by this perturbation. The three exceptional cases can be treated even more easily using the computer to have an explicit expression for $U$ and its eigenvalues at any point of the manifold.

\vspace{0.5cm}

Recall that the orbifold cohomology of $\mathbb{P}^1_{p,q,r}$ as a vector space is just the singular homology of the inertia orbifold of $\mathbb{P}^1_{p,q,r}$, i.e. topologically the disjoint union of a sphere $S^2$ and $(p+q+r-3)$ isolated points which we denote $X_j,Y_k,Z_l$ with  $j=1,\ldots,p-1$, $k=1,\ldots,q-1$, $l=1,\ldots,r-1$. We choose the homogeneous basis of cohomology classes for $QH^*_\mathrm{orb}(\mathbb{P}^1_{2,2,r})$ given by
$$1 \in H^0(S^2)\hspace{0.5cm},\hspace{0.5cm} p \in H^2(S^2)$$
$$x_1\in H^0(X_1)\hspace{0.5cm},\hspace{0.5cm} y_1 \in H^0(Y_1)$$
$$z_1\in H^0(Z_1),\ \ldots,\ z_{r-1}\in H^0(Z_{r-1})$$
Let us first compute the quantum algebra structure on $QH^*_\mathrm{orb}(\mathbb{P}^1_{2,2,r})$ at a point $m$ whose flat coordinates (components on the above basis of cohomology classes) are zero, with the exception of the components $a$ and $b$ along $x_1$ and $y_1$, which are left generic.

This structure is a deformation of the ordinary orbifold cup product on $H^*_\mathrm{orb}(\mathbb{P}^1_{p,q,r})$, which is given by
$$x_1\cdot y_1=x_1 \cdot z_k=y_1\cdot z_k=0\ \mathrm{if}\ k=1,\ldots,r-1$$
$$z_{k_1}\cdot z_{k_2}=x_{k_1+k_2}\ \mathrm{if}\ k_1+k_2\leq r-1$$
$$2x_1^2=2y_1^2=rz_1^r=p$$
It follows that, as rings,
$$H^*_{orb}(\mathbb{P}^1_{2,2,r})\simeq \C [x,y,z]/(xy,xz,yz,2x^2-2y^2,2x^2-rz^r,2y^2-rz^r)$$
if we identify $x_1=x$, $y_1=y$, $z_1=z$.

The main result for determining the algebra structure on $T_mQH^*_\mathrm{orb}(\mathbb{P}^1_{2,2,r})= H^*_\mathrm{orb}(\mathbb{P}^1_{2,2,r})\otimes \C[q]$ at $m$ is the following.
\begin{lemma}
$$x_1\bullet_m z_1=2qy_1+bq$$
$$y_1\bullet_m z_1=2qx_1+aq$$
$$x_1\bullet_m y_1=r(qz_{r-1}+q^3z_{r-3}+q^5z_{r-5}+\ldots+q^{2\lfloor \frac{r+1}{2}\rfloor-1}z_{r-(2\lfloor \frac{r+1}{2}\rfloor-1)})$$
$$z_1\bullet_m z_1=z_2+2q^2$$
$$z_1\bullet_m z_2=z_3+q^2z_1$$
$$\ldots$$
$$z_1\bullet_m z_{r-2}=z_{r-1}+q^2z_{r-3}$$
$$z_1\bullet_m z_{r-1}=p+rq^2z_{r-2}$$
\end{lemma}
\begin{proof}
Recall firstly that the quantum product agrees with the grading of orbifold cohomology given by
$$\mathrm{deg}\, z_i=\frac{2i}{r}\hspace{1cm}\mathrm{deg}\, x_1=\mathrm{deg}\, y_1=1 \hspace{1cm} \mathrm{deg}\, p=2$$
$$\mathrm{deg}\, q=\frac{2}{r} \hspace{1cm} \mathrm{deg}\, a=\mathrm{deg}\, b=1$$
This gives selection rules on the form of the product. For instance the most general possibility for the product $x_1\bullet_m z_1$ is $$x_1\bullet_m z_1=2c_1 x_1q+2c_2 y_1q+rc_3 a z_1+rc_4 b z_1+c_5 a q + c_6 b q+r\sum_ic_{7,i}z_i q^{\frac{r+2}{2}-i}$$ where the last sum is only present if $r$ is even and, by the explicit form of the Poincar\'e pairing and by definition of quantum multiplication (with the usual Gromov-Witten invariants correlator bracket notation $<\ldots>_{g,k,d}$, where the indices are respectively genus $g$, number of marked points $k$ and degree $d$),
$$c_1=<x_1,z_1,x_1>_{0,3,1}=0$$
$$c_2=<x_1,z_1,y_1>_{0,3,1}=1$$
$$c_3=<x_1,z_1,z_{r-1},x_1>_{0,4,1}=0$$
$$c_4=<x_1,z_1,z_{r-1},y_1>_{0,4,1}=0$$
$$c_5=<x_1,z_1,p,x_1>_{0,4,1}=0$$
$$c_6=<x_1,z_1,p,y_1>_{0,4,1}=1$$
$$c_{7,i}=<x_1,z_1,z_i>_{0,4,\frac{r+2}{2}-i}=0$$
The vanishing of $c_1$, $c_3$, $c_4$, $c_5$ and the $c_{7,i}$ is proven (see also \cite{MT}) by considering a stable map $f$ in the moduli space relevant for the corresponding correlator and the pull back via $f$ of the three line bundles $L_x$, $L_y$, $L_z$ which generate the Picard group of $\mathbb{P}^1_{2,2,r}$. The holomorphic Euler characteristics $\chi$ of these bundles have to be integer and can be computed explicitly via Riemann-Roch. For instance, for a map $f$ in the moduli space relevant for the correlator $c_4$, one has
$$\chi(f^*L_x)=1+\frac{1}{2}-\frac{1}{2}$$
$$\chi(f^*L_y)=1+\frac{1}{2}-\frac{1}{2}$$
$$\chi(f^*L_z)=1+\frac{1}{r}-\frac{1}{r}-\frac{r-1}{r}$$
and since the last one is not an integer we deduce that the moduli space is empty and the correlator vanishes. For the same reason $c_1$, $c_3$,$c_5$ and the $c_{7,i}$ vanish, while the correlators $<x_1,z_1,y_1>_{0,3,1}$ and $<x_1,z_1,p,y_1>_{0,4,1}$ are clearly $1$.

The same technique can be used for all the other relations of the statement but the correlators $<z_1,z_i,z_{i-d+1}>_{0,3,d}$ appearing in $z_1\bullet_m z_i$. In this case the Euler characteristic only selects $d$ to be even and in order to kill all the terms with $d\neq 2$ one needs to consider a bit more carefully the maps in the relevant moduli space. The key observation here is that any constant component of the map needs extra marked points to be stabilized, while non-constant components need extra marked orbifold points every time they induce a branched covering of $S^2$ which has local degree not multiple of $2,2,r$ locally over the orbifold points of stabilizer $\mathbb{Z}_2,\mathbb{Z}_2,\mathbb{Z}_r$. Plugging this information into Riemann-Hurwitz relation in the same spirit as in the proof of Theorem \ref{polyP1}, one gets the desired vanishing result.
\end{proof}

The above lemma implies
$$T_mQH^*_\mathrm{orb}(\mathbb{P}^1_{2,2,r})\simeq \frac{\C[[q]] [x,y,z]}{ \begin{pmatrix} xy-rqz^{r-1}+r\sum_{k=1}^{\lfloor\frac{r-1}{2}\rfloor}(-1)^{k-1}(r-2k)\frac{(r-k-1)!}{k!(r-2k)!}q^{2k+1}z^{r-2k-1} \\ xz-2qy-bq \\ yz-2qx-aq \end{pmatrix}}$$
where $x_1=x$, $y_1=y$ and $z_1=z$ and where we can set $q$ to any nonzero complex number.

Using the flat coordinates of Lemma \ref{flat} it is easy to see that, at our point $m$, $c_{r-2k}=(-1)^{k}r \frac{(r-k-1)!}{k!(r-2k)!}$ for $k=1,\ldots\lfloor\frac{r-1}{2}\rfloor$, which implies that, for q=1 and for every $a$ and $b$, the one above is precisely the algebra $T_mM_{2,2,r}\simeq\C[x,y,z]/J_{F_m}$ where $$F_m(x,y,z)=-xyz+x^2+ax+y^2+by+z^r+r\sum_{k=1}^{\lfloor\frac{r-1}{2}\rfloor} (-1)^{k}\frac{(r-k-1)!}{k!(r-2k)!}q^{2k+1}z^{r-2k}$$

This way we can directly compute the linear operator $U$ of (quantum) multiplication by $E$ at $m$, whose matrix with respect to the flat coordinates we still denote (with a little abuse of notation) $U=U_0+V$, where $V(a=0,b=0)=0$ and
\begin{equation*}V=\left(
\begin{array}{ccccccccccc}
0 &                &                &             &        &        &          &           &          &        &  A \\
  & -\frac{a^2}{4} &      0         &             &        &        &          &           &          &        &    \\
  &       0        & -\frac{b^2}{4} &             &        &        &          &           &          &        &  \\
  &                &                &      0      &        &        &          &           &          &   B_{\lfloor\frac{r-1}{2}\rfloor}  &  \\
  &                &                &    \vdots   & \ddots &        &          &           &  \ldots  &        &  \\
  &                &                &      0      & \ldots &   0    &          &   B_1     &          &        &  \\
  &                &                &      0      & \ldots & \ldots &     0    &           &          &        &  \\
  &                &                &    \vdots   & \ldots & \ldots &   \vdots &     0     &          &        &  \\
  &                &                &    \vdots   & \ldots & \ldots &   \vdots &  \vdots   &  \ddots  &        &  \\
  &                &                &      0      & \ldots & \ldots &     0    &     0     &  \ldots  &    0   &  \\
0 &                &                &             &        &        &          &           &          &        & 0\\
\end{array}\right)
\end{equation*}
for $r$ even and
\begin{equation*}V=
\begin{pmatrix}
0 &                &                &             &           &          &           &          &        &  A \\
  & -\frac{a^2}{4} &      0         &             &           &          &           &          &        &    \\
  &       0        & -\frac{b^2}{4} &             &           &          &           &          &        &  \\
  &                &                &      0      &           &          &           &          &   B_{\lfloor\frac{r-1}{2}\rfloor}  &  \\
  &                &                &      0      &  \ddots   &          &           &  \ldots  &        &  \\
  &                &                &      0      &  \ldots   &     0    &   B_1     &          &        &  \\
  &                &                &      0      &  \ldots   &     0    &    0      &          &        &  \\
  &                &                &  \vdots     &  \vdots   &  \vdots  &   \vdots  &  \ddots  &        &  \\
  &                &                &      0      &  \ldots   &     0    &    0      &  \ldots  &    0   &  \\
0 &                &                &             &           &          &           &          &        & 0\\
\end{pmatrix}
\end{equation*}
for $r$ odd, with $A=2r\frac{r}{2}(a^2+b^2)$, $B_k= \frac{r+2k}{r}k(a^2+b^2)$ for $r$ even and $A=rab$, $B_k= \frac{r+2k-1}{r}(2k-1)ab$ for $r$ odd, and we won't need to know the non-specified entries. While

\begin{equation*}U_0=
\begin{pmatrix}
0 &                &                &             &           &          &                     &        &  4r \\
  &       2        &       0        &             &           &          &                     &        &    \\
  &       0        &       2        &             &           &          &                     &        &  \\
  &                &                &             &           &          &                     &   2    &  \\
  &                &                &             &           &          &             \ldots  &        &  \\
  &                &                &             &           &          &                     &        &  \\
  &                &                &             &           &     2    &                     &        &  \\
  &                &                &             &  \ldots   &          &                     &        &  \\
  &                &                &      2      &           &          &                     &        &  \\
 \frac{1}{r} &     &                &             &           &          &                     &        & 0\\
\end{pmatrix}
\end{equation*}
for $r$ even and
\begin{equation*}U_0=
\begin{pmatrix}
0 &                &                &             &           &          &           &          &        &  4r \\
  &       0        &       2        &             &           &          &           &          &        &    \\
  &       2        &       0        &             &           &          &           &          &        &  \\
  &                &                &             &           &          &           &          &   2    &  \\
  &                &                &             &           &          &           &  \ldots  &        &  \\
  &                &                &             &           &          &   2       &          &        &  \\
  &                &                &             &           &     2    &           &          &        &  \\
  &                &                &             &  \ldots   &          &           &          &        &  \\
  &                &                &      2      &           &          &           &          &        &  \\
 \frac{1}{r} &     &                &             &           &          &           &          &        & 0\\
\end{pmatrix}
\end{equation*}
for $r$ odd, but here the non-specified entries are zeros.

Now we use perturbation theory to show that $U$ has distinct eigenvalues. In fact let $v_1,\ldots,v_{p+q+r-1}$ be a basis of eigenvectors of $U_0$, and $\tilde{V}$ the matrix of $V$ in this basis. If $(U_0+\epsilon V) v_i = \lambda(\epsilon) v_i$, then $\lambda(\epsilon)=\lambda(0)+\epsilon \tilde{V}_{ii}+O(\epsilon^2)$. Checking that, for $a\neq b$, the perturbation completely removes the degeneracy of eigenvalues, using the above explicit expressions, is straightforward. This completes the proof of Theorem \ref{mirror}.

\section{Symplectic Field Theory of Seifert fibrations}\label{SFT}

What we are going to explain now is a generalization of Proposition $2.9.2$ of \cite{EGH} to the case of oriented Seifert fibrations over orbifold Riemann surfaces. For simplicity we will consider just the case of an oriented Seifert $S^1$-orbibundle $\pi:V\to P$ over a $\mathbb{P}^1$-orbifold $P=\mathbb{P}^1_{(z_1,\alpha_1),\ldots,(z_a,\alpha_a)}$ with uniformizing systems at $z_i$ given by $z\mapsto z^{\alpha_i}$, with Seifert invariant $(c,\beta_1,\ldots,\beta_a)$ (here $c$ is the first Chern class of the orbibundle $V$ and $b=c-\sum\frac{\beta_i}{\alpha_i}$ is the first Chern class of its de-singularization $|V|$, as in \cite{CR}). This is actually a sort of mildly singular Hamiltonian structure of fibration type, since the manifold $V$ itself is still smooth, but the Reeb orbit space is a symplectic orbifold with $\mathbb{Z}_n$-singularities. However it satisfies the \emph{Morse-Bott} condition of \cite{B}. Moreover we stick to the rational ($g=0$) SFT. Here we use Chen and Ruan's theory and terminology of orbifolds \cite{CR},\cite{CR1}.

The cylindrical cobordism $V\times\mathbb{R}$ can be seen as the total space of the complex line orbibundle $L$ associated with $V\to P$, which is actually holomorphic for a proper choice of the almost complex structure $J$, and with the zero section (containing the singular locus of $L$) removed. Following \cite{EGH} we want to project a SFT-curve $u:\mathbb{P}^1-\{x_1,\ldots,x_s\}\to L$ along $\pi$ to obtain a stable orbifold GW-curve $\bar{u}$ to $P$. This is trivial outside $\pi^{-1}(z_i)$, while the following Lemma ensures that the domain curve $\mathbb{P}^1-\{x_1,\ldots,x_s\}$ can be given in a unique way an orbifold structure at the punctures $\{x_1,\ldots,x_s\}$ in such a way that we get a genuine stable GW-map to $P$.
\begin{lemma}
Let $u:\Sigma\to\Sigma'_{(z_1,\alpha_1),\ldots,(z_a,\alpha_a)}$ be a non-constant holomorphic map between a Riemann surface and a complex orbicurve. Then there is a unique orbifold structure on $\Sigma$ (which we denote $\Sigma_u$) and a unique germ of $C^\infty$-lift $\tilde{u}$ of $u$ (which is regular) such that the group homomorphism at each point is injective.
\end{lemma}
\begin{proof}
There are finitely many $x_i\in\Sigma$ such that $u(x_i)=z_i$. An orbifold structure at $x_i$ and the corresponding $C^\infty$-lift of $u$ fit into the diagram
\begin{equation*}
\xymatrix{D       \ar[r]^{z^q} \ar[d]_{z^p}& D       \ar[d]^{z^{\alpha_i}} \\
          D_{x_i} \ar[r]^{z^n}             & D_{z_i} }
\end{equation*}
where $n$ is assigned with the map $u$. So we have $pn=q\alpha_i$. Now, if $p$ and $q$ have some common factor, then the group homomorphism is not injective (the common factor $k\in\mathbb{Z}_p$ is sent to $0\in\mathbb{Z}_{\alpha_i}$), then $p$ and $q$ are uniquely determined by $n$ and $\alpha_i$ thanks to the injectivity condition. The map is obviously regular by construction.
\end{proof}

This way we get a fibration between the relevant moduli spaces of holomorphic maps
$$\mathrm{pr}:\mathcal{M}^A_{0,r,s}(L)\to\mathcal{M}^A_{0,r+s}(P)$$
where $A\in H_2(P,\mathbb{Z})=\mathbb{Z}$.
Using the orbicurve version of the correspondence between effective divisors and line bundles on $P$, one is able to identify the fiber $\mathrm{pr}^{-1}(\bar{u})$. As in the smooth case, $u$ can be reconstructed from $\bar{u}$ by assigning a meromorphic section of $\bar{u}^*(L)$ which is determined, up to an $S^1$-symmetry, by its divisor of poles and zeros, hence by $s$ of the $r+s$ marked points on $\mathbb{P}^1_{\bar{u}}$ and a sequence of integers $(k_1,\ldots,k_s)$ such that $$\sum_{i=1}^s\frac{k_i}{m_i}=c_1(\bar{u}^*(L))=\mathrm{deg}(\bar{u})\, c_1(L)=\mathrm{deg}(\bar{u})\, \left( c_1(|L|)+\sum_{j=1}^{a}\frac{\beta_i}{\alpha_i}\right),$$ where $m_i$ is the orbifold multiplicity of the $i$-th of the $s$ points of the divisor on $\mathbb{P}^1_{\bar{u}}$ and $|L|$ is the de-singularization of $L$ (see e.g. \cite{CR}).

Let now $\Delta_1,\ldots,\Delta_b$ be a basis of $H^*_\mathrm{orb}(P)$ such that the system of forms $\tilde{\Delta}_j:=\pi^*(\Delta_j)$, $j=1,\ldots,c<b$ generate $\pi^*(H^*_{\mathrm{orb}}(P))=\pi^*(H^*(P))\subset H^*(V)$, and the forms $\tilde{\Theta}_1,\ldots,\tilde{\Theta}_d$ complete it to a basis of $H^*(V)$. Notice also that $H_2(V,\mathbb{Z})\cong H^1(V,\mathbb{Z})=0$ (if $V$ is not trivial).\\
Extending the proof in \cite{B}\cite{EGH} to obtain the following generalization is now a trivial matter.
\begin{prop}
Let $\mathbf{f}_P(\sum t_i\Delta_i,z)$ be the genus $0$ orbifold Gromov-Witten potential of $P$ and $\mathbf{h}_V(\sum t_i\tilde{\Delta}_i+\sum \tau_j\tilde{\Theta}_j,q,p)$ the rational SFT potential of $V$ (as a framed Hamiltonian structure of fibration type). Let
$$\mathbf{h}^j_V(t,q,p)=\left.\frac{\partial \mathbf{h}_V}{\partial \tau_j}\left(\sum_{1}^{c} t_i\tilde{\Delta}_i+\tau_j\tilde{\Theta}_j,q,p\right)\right|_{\tau_j=0}$$
$$\mathbf{f}^j_P(t;z)=\left.\frac{\partial \mathbf{f}_P}{\partial s}\left(\sum_{1}^{b} t_i\Delta_i+s\pi_*\tilde{\Theta}_j,z\right)\right|_{s=0}$$
for $j=1,\ldots,d$ and where $\pi_*$ denotes integration along the fibers of $V$. Then we have
$$\mathbf{h}^j_V(t,q,p)=\frac{1}{2\pi\alpha_1\ldots\alpha_a}\int_{0}^{2\pi\alpha_1\ldots\alpha_a}\mathbf{f}^j_P(t_1+u_1(x),\ldots,t_b+u_b(x),u_{b+1}(x),\ldots,u_c(x);e^{-\mathrm{i}c_1(V)x})\di x$$
where $$u_n(x):=\sum_{k=1}^\infty \left(q_{k\iota_n,n}\, e^{-\mathrm{i}\,k\iota_n x}+ p_{k\iota_n ,n}\, e^{\mathrm{i}\,k\iota_n x}\right)\qquad n=1,\ldots,b$$
and $\iota_n$ denotes the degree shifting (in the sense of orbifold cohomology, see \cite{CR}\cite{CR1}) of the twisted sector of $P$ where the form $\Delta_n$ is defined.
\end{prop}

More explicitly, up to degree shifting and just as a vector space, the orbifold cohomology of $P$ is the singular cohomology (with $\mathbb{R}$ coefficients) of the inertia orbifold of $P$, hence just the disjoint union of a sphere $S^2$ and $[(\alpha_1-1)+\ldots +(\alpha_a-1)]$ points, labeled by a couple of indices $(i=1,\ldots,a;l=1,\ldots,\alpha_i-1)$. Let us choose the $(\alpha_1+\ldots+\alpha_a-a+2)$ elements of the basis of $H^*_\mathrm{orb}(P)$ in the following way: $\Delta_1=1,\Delta_2=[\omega]\in H^*(S^2)$ and $\Delta_{(i,l)}=1\in H^0(\mathrm{pt}_{(i,l)})$. Then $\iota_1=\iota_2=0$ and $\iota_{(i,l)}=\frac{l}{\alpha_i}$

Of course, as in the smooth case of \cite{EGH}, the Poisson structure and grading on the graded Poisson algebra $\mathfrak{U}$, relevant for rational SFT and where $\mathbf{h}_V$ is defined, is given again and in the same way by the orbifold Poincar\'e pairing and grading in $H^*_{\mathrm{orb}}(P)$ (see \cite{CR}).

\vspace{0.5cm}

Thanks to this generalization we can use the Gromov-Witten potentials of $\mathbb{P}^1$-orbifolds computed in the previous section to obtain explicit expressions for the SFT-Hamiltonians of some interesting contact $3$-manifolds, which are Seifert fibrations over these orbifolds. The $A$-case gives the general Lens spaces, the $D$-cases are Prism manifolds, while among the exceptional $E$-cases one can find quotients of the Poincar\'e sphere and more exotic manifolds. We plan to investigate these examples in subsequent publications.

\end{document}